\documentclass[journal,letterpaper]{IEEEtran}
\usepackage{amsmath,amsfonts,amssymb,epsfig,url}
\newenvironment{myproof}[1][Proof]{{\em #1:} }{\hspace*{\fill} \rule{0.5em}{0.5em}} 
\newcommand{\Reals}{\mathbb{R}}
\newcommand{\NN}{\mathbb N}
\newtheorem{theorem}{\textbf{Theorem}}
 
\newtheorem{definition}{\textbf{Definition}}

\newcommand{\mybar}[1]{\overline{#1}}
\newcommand{\beql}[1]{\begin{equation}\label{#1}}
\newcommand{\eeq}{\end{equation}}

\renewcommand{\vec}[1]{\mbox{\boldmath$#1$}}

\begin{document}

\title{Constant Weight Codes: A Geometric Approach Based on Dissections}

\author{Chao Tian,~\IEEEmembership{Member,~IEEE,}
Vinay A. Vaishampayan~\IEEEmembership{Senior Member,~IEEE}
and
N. J. A. Sloane~\IEEEmembership{Fellow,~IEEE}
\thanks{Chao Tian is with EPFL Lausanne. Vinay A. Vaishampayan and N. J. A. Sloane are with AT\&T Shannon Laboratory, Florham Park, NJ. This work was done while Chao Tian was visiting AT\&T Shannon Laboratory.}}
     \maketitle

\begin{abstract}
We present a novel technique for encoding and decoding
constant weight binary codes that uses a
geometric interpretation of the codebook. 
Our technique is based on embedding the codebook in a Euclidean space 
of  dimension equal to the weight of the code.
The encoder and decoder mappings  are
then interpreted as a bijection between a certain
hyper-rectangle and a polytope
in this Euclidean space. An inductive dissection algorithm
is developed for constructing such a bijection.  
We prove that the algorithm is correct and then analyze
its complexity. The complexity depends on the weight of the 
code, rather than on the block length as in other algorithms.
This approach is advantageous when the weight is smaller than the square root of the block length.

\end{abstract}

\section{Introduction}
\label{sec:intro}
We consider the problem of encoding and decoding binary codes
of constant Hamming weight $w$ and block length $n$.
Such codes are useful in a  variety of applications: a few examples are fault-tolerant circuit design and computing~\cite{PrSt1980},
pattern generation for circuit testing~\cite{TaWo1983}, identification coding~\cite{VeWe1993} and optical overlay networks~\cite{VaiFeu:2005}.

The problem of interest is that of designing the encoder and decoder,
i.e., the problem of mapping all binary (information) vectors of a given length onto a subset of length-$n$ vectors of constant Hamming weight $w$ 
in a one-to-one manner. 
In this work, we propose a novel geometric method in which information and
code vectors are represented by vectors in $w$-dimensional Euclidean space, 
covering polytopes for the two sets are identified, and a one-to-one mapping 
is established by dissecting the  covering polytopes in a specific manner.
This approach results in an invertible integer-to-integer mapping, 
thereby ensuring unique decodability. 
The proposed algorithm has a natural recursive structure, 
and an inductive proof
is given for unique decodability. The issue of efficient  encoding and decoding is also addressed. We show that
 the proposed algorithm has complexity $O(w^2)$,
 where $w$ is the weight of the codeword, independent of the codeword length.

Dissections are of considerable interest in geometry,
partly as a source of puzzles, but more importantly because 
they are intrinsic to the notion of volume.
Of the $23$ problems posed by David Hilbert
at the International Congress of Mathematicians in 1900,
the third problem dealt with dissections.
Hilbert asked for a proof that there are two
tetrahedra of the same volume with the property
that it is impossible to dissect one into a finite
number of pieces that can be rearranged to give the other,
i.e., that the two tetrahedra are not equidecomposable.
The problem was immediately solved by Dehn~\cite{Dehn1900}.
In 1965, after $20$ years of effort, Sydler~\cite{Sydl1965}
completed Dehn's work. 
The Dehn-Sydler theorem states that a necessary and sufficient condition
for two polyhedra to be equidecomposable is that they have the same
volume and the same Dehn invariant.
This invariant is a certain function of the edge-lengths and dihedral angles of the polyhedron.
An analogous theorem holds in four dimensions (Jessen~\cite{Jess1968}),
but in higher dimensions it is known only that equality of the Dehn invariants
is a necessary condition. In two dimensions any two polygons of 
equal area are equidecomposable, a result due to
Bolyai and Gerwein (see Boltianskii~\cite{Bolt1978}).
Among other books dealing with the classical dissection problem
in two and three dimensions we mention in
particular Frederickson~\cite{Fred1997},
Lindgren~\cite{Lind64}
and Sah~\cite{Sah1979}.

The remainder of the paper is organized as follows. 
We provide background and review relevant previous work
in Section~\ref{sec:review}. 
Section \ref{sec:geo} describes our
geometric approach and gives some low-dimensional examples.
Encoding and decoding algorithms are then given in Section~\ref{sec:alg},
and the correctness of the algorithms is established.
Section \ref{sec:con} summarizes the paper.

\section{Background and Previous Methods}
\label{sec:review}

Let us denote the Hamming weight of a length-$n$ binary sequence
$\vec{s} := (s_1,s_2,\ldots,s_n)$ by $w(s) := |\{s_i:s_i=1\}|$,
where $|\cdot|$ is the cardinality of a set. 

\begin{definition}
An $(n,w)$ constant weight binary code ${\cal C}$ is a set of length-$n$ sequences
such that any sequence $\vec{s}\in {\cal C}$ has weight $w(\vec{s})=w$.
\end{definition}

If  $\cal{C}$ is an $(n,w)$ constant weight code,
then its rate $R:=(1/n)\log_2|{\cal C}|\leq R(n,w):=(1/n) \log_2 \binom{n}{w}$.
For fixed $\beta$ and $w=\left\lfloor \beta n\right \rfloor$, we have
\begin{equation}
\mybar{R}:=\lim_{n\rightarrow \infty} R(n,w) = h(\beta) \,,
\end{equation}
where $h(\beta) := -\beta\log_2(\beta)-(1-\beta)\log_2(1-\beta)$ is the entropy function.
Thus $\mybar{R}$ is maximized when $\beta=1/2$, i.e., the asymptotic rate is highest when the code is balanced.

The (asymptotic) efficiency of a code  relative
to an infinite-length code with 
the same weight to length ratio $w/n$, given by
$\eta := R/\mybar{R}$, can be written
as $\eta=\eta_1 \mybar{\eta}$ where  $\eta_1:= R/R(n,w)$
and $\mybar{\eta}:= R(n,w)/\mybar{R}$. 
The first term, $\eta_1$, is the efficiency of a particular code
relative to the best possible code with the same length and weight;
the second term, $\mybar{\eta}$, is the efficiency of the best
finite-length code relative to the best infinite-length code. 

\begin{figure}[t]
\begin{center}
\includegraphics[width=7.5cm]{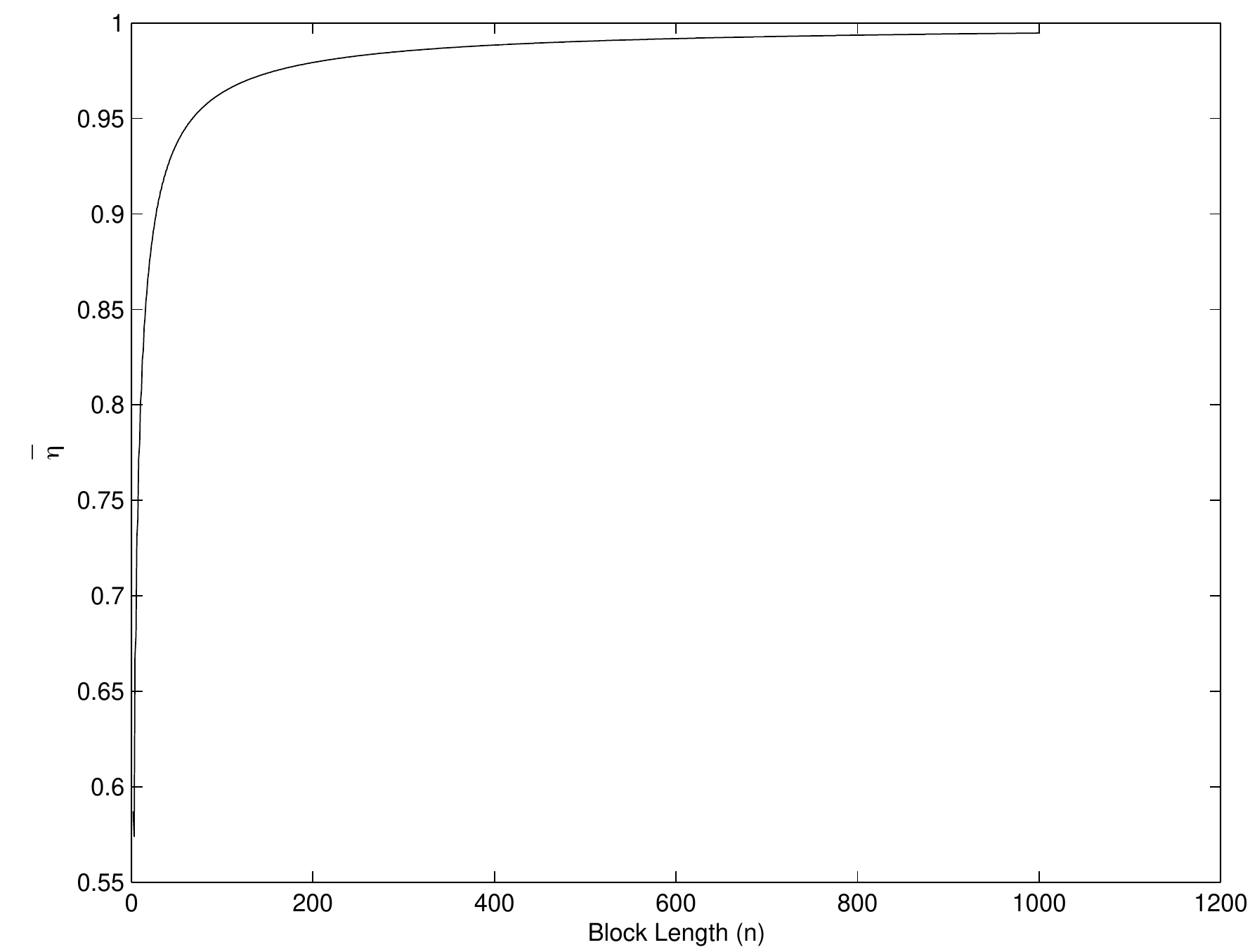}
\end{center}
\caption{Efficiency $\mybar{\eta}$ as a function of block length when
 $\beta=1/2$}
\label{fig:efficiency}
\label{fig1}
\end{figure}

 From Stirling's formula we have 
\begin{equation}
\mybar{\eta}\approx 
1 - \frac{ \log_2 (2 \pi n\beta (1-\beta))}{2 n h(\beta)}.
\end{equation}
A plot of $\mybar{\eta}$ as a function of $n$ is given in Fig.~\ref{fig:efficiency} for
$\beta=1/2$. 
The slow convergence visible here
is the reason one needs
codes with large block lengths.

Comprehensive tables and construction techniques for binary constant weight
codes can be found in~\cite{BSSS1990} and the references therein.
However, the problem of finding efficient encoding and decoding algorithms
has received considerably less attention. 
We briefly discuss two previous  methods that are relevant to our work. 
The first, a general purpose technique based on the idea of
lexicographic ordering and enumeration of codewords in a codebook
(Schalkwijk~\cite{Schl1972}, Cover~\cite{Cove1973})
is an example of ranking/unranking algorithms that are well studied
in the combinatorial literature
(Nijenhuis and Wilf~\cite{NW1978}).
We refer to this as the {\em enumerative} approach.
The second (Knuth~\cite{Knut1986}) is a special-purpose, highly efficient technique
that works for balanced codes, i.e.,
when $w=\lfloor(n/2)\rfloor$, and is referred to as the {\em complementation} method.

The enumerative approach orders the codewords lexicographically
(with respect to the partial order defined by $0 < 1$),
as in a dictionary. The encoder computes the codeword from its dictionary index,
and the decoder computes the dictionary index from the codeword. 
The method is effective because there is a simple formula
involving binomial coefficients for computing the lexicographic index of a codeword.
The resulting code is fully efficient in the sense that $\eta_1=1$.
However, this method requires the computation of the exact values of
binomial coefficients $\binom{n}{k}$, and requires registers
of length $O(n)$, which limits its usefulness.

An alternative is to use arithmetic
coding (Rissanen and Langdon~\cite{RiLa1979}, Rissanen~\cite{Riss1979};
see also Cover and Thomas~\cite[\S 13.3]{CoTh2006}).
Arithmetic coding is an efficient variable length source coding 
technique for finite alphabet sources. 
Given a source alphabet and a simple probability 
model for sequences $\vec{x}$,
let $p(\vec{x})$ and $F(\vec{x})$ denote the probability distribution and
cumulative distribution function, respectively. 
An arithmetic encoder represents $\vec{x}$ 
by a number in the interval $(F(\vec{x})-p(\vec{x}),F(\vec{x})]$.  
The implementation of such a coder
can also run into problems with very long registers,
but elegant finite-length implementations are known
and are widely used
(Witten, Neal and Cleary~~\cite{WNC1987}).
For constant weight codes,
the idea is to reverse the roles of encoder and decoder,
i.e., to use an arithmetic decoder as an encoder and an
arithmetic encoder as a constant weight decoder
(Ramabadran~\cite{Rama1990}). 
Ramabadran gives an efficient  algorithm based on an adaptive
probability model,
in the sense that the probability that the incoming  bit is a $1$
depends on the number of $1$'s
that have already occurred. 
This approach successfully overcomes the
 finite-register-length constraints associated with computing the binomial
coefficients and  
the resulting efficiency is often very high,
 in many cases the loss of information being at most one bit.
The encoding complexity of the method is $O(n)$.

Knuth's complementation method~\cite{Knut1986}  relies on the key observation that
if the bits of a  length-$n$ binary
sequence are complemented sequentially,
 starting from the beginning,
 there must be a point
 at which the weight is equal to $\lfloor{n/2}\rfloor$. 
Given the transformed sequence,
 it is possible to recover
the original sequence by specifying how many bits were complemented 
(or the weight of  the original sequence). 
This information is provided by a (relatively  short)
constant weight check string,
and the resulting code consists of the transformed sequence followed by
the constant weight check bits.  In a series of papers,
 Bose and colleagues extended Knuth's method in various ways,
 and determined the limits of this approach
(see \cite{YoBo2003} and references therein).
The method is simple and efficient,
 and even though the overall complexity is $O(n)$,
for $n=100$
 we found it to be eight times as fast
as the method based on arithmetic codes. However,
 the method only works for balanced codes,
which restricts its applicability.

The two techniques that we have described above both have complexity 
that depends on the length $n$ of the codewords.
In contrast, the complexity of our algorithm depends
only on the weight $w$, which makes it more  suitable for codes with relatively low weight. 

As a final piece of background information, we define
what we mean by a dissection. 
We assume the reader is familiar with the terminology of polytopes
(see for example
Coxeter~\cite{Coxe1973},
Gr\"unbaum~\cite{Grun2003},
Ziegler~\cite{Zieg1995}).
Two polytopes $P$ and $Q$ in $\Reals^w$ are  said to be
\emph{congruent} if $Q$ can be obtained from $P$ by a translation, 
a rotation and possibly a reflection in a hyperplane.
Two polytopes $P$ and $Q$ in $\Reals^w$ are  said to be
\emph{equidecomposable} if they can be decomposed into 
finite sets of polytopes $P_1, \ldots, P_t$ and
$Q_1, \ldots, Q_t$ , respectively, for some positive
integer $t$, such that $P_i$ and $Q_i$ are
congruent for all $i=1, \ldots,t$
(see Frederickson~\cite{Fred1997}).
That is, $P$ is the disjoint union of the
polytopes $P_i$, and similarly for $Q$.
If this is the case then we say that $P$ can be 
\emph{dissected} to give $Q$ (and that $Q$ can be
dissected to give $P$).

Note that we allow reflections in the dissection:
there are at least four reasons 
for doing so. (i) It makes no difference to
the \emph{existence} of the dissection, since
if two polytopes are equidecomposable using reflections
they are also equidecomposable without using
reflections. This is a classical theorem in two and three
dimensions \cite[Chap.~20]{Fred1997} and the proof is easily
generalized to higher dimensions.
(ii) When studying congruences, it is simpler not to have to worry about
whether the determinant of the orthogonal matrix has determinant
$+1$ or $-1$. 
(iii) Allowing reflections often reduces the number of pieces.
(iv) Since our dissections are mostly in dimensions
greater than three, the question of ``physical realizability''
is usually irrelevant.
Note also that we do not require that the $P_i$ can be obtained
from $P$ by a succession of cuts along infinite hyperplanes.
All we require is that $P$ be a disjoint union of the $P_i$.

One final technical point: when defining dissections
using coordinates, as in Eqns.
(\ref{eq2Da}), (\ref{eq2Db}) below,
we use a mixture of $\leq$
and $<$ signs in order to have unambiguously defined maps.
This is essential for our application.
On the other hand, it means that the ``pieces''
in the dissection may be missing certain boundaries.
It should therefore be understood that if we were
focusing on the dissections themselves,
we would replace each piece by its
topological closure.

For further information about dissections see
the books mentioned in Section \ref{sec:intro}.

\section{The Geometric Interpretation}
\label{sec:geo}

In this section, we first consider the problem of
encoding and decoding a binary constant weight code of weight $w=2$ 
and arbitrary  length $n$, i.e., where there are only two bits set
to $1$ in any codeword. 
Our approach is based on the fact that vectors
of weight two can be represented as points
in two-dimensional Euclidean space, and can 
be scaled, or normalized, to lie in a  right triangle.
This approach is then  extended, first to weight $w=3$, and then
to arbitrary weights $w$. 

For any weight $w$ and block length $n$,
let ${\cal C}_w$
denote the set of all weight $w$ vectors, with
$|{\cal C}_w| = \binom{n}{w}$.
Our codebook ${\cal C}$ will be a subset of ${\cal C}_w$,
and will be equal to ${\cal C}_w$ for a fully efficient code, 
i.e., when $\eta_1=1$.
We will represent a codeword
by the $w$-tuple $\vec{y'}  :=  (y'_1,y'_2,\ldots,y'_w)$,
$1 \leq y'_1 < y'_2 < \ldots < y'_w \leq n$, 
where $y'_i$ is the position of the $i$th $1$ in the codeword,
counting from the left. 
If we normalize these indices $y'_i$ by dividing them by $n$,
the codebook ${\cal C}$ becomes a discrete subset of  
the polytope $T_w$, the convex hull of the points
$ 0^w, 0^{w-1}1, 0^{w-2}11, \ldots, 01^{w-1}, 1^w$.
$T_2$ is a right triangle,
$T_3$ is a right tetrahedron and
in general we will call $T_w$ 
a {\em unit orthoscheme}\footnote{An
{\em orthoscheme} is a $w$-dimensional simplex
having an edge path consisting of $w$ totally 
orthogonal vectors (Coxeter~\cite{Coxe1973}).
In a {\em unit orthoscheme} these edges
all have length $1$.}.

The set of inputs to the encoder will be denoted by ${\cal R}_w$:
we assume that
this consists of $w$-tuples $\vec{y}  :=  (y_1,y_2,\ldots,y_w)$
which range over a $w$-dimensional hyper-rectangle or ``brick''.
After normalization by dividing the $y_i$ by $n$, we may assume that
the input vector is a point in the hyper-rectangle or ``brick''
$$
B_w:=[0,1) \times [1-1/2,1) \times \ldots \times [1-1/w,1) \,.
$$
We will use $\vec{x} := (x_1,x_2,\ldots,x_w) = \vec{y}/n \in B_w$
and $\vec{x'} := (x'_1,x'_2,\ldots,x'_w) = \vec{y'}/n \in T_w$
to denote the normalized versions
of the input vector and codeword, respectively, defined by
$x_i := y_i/n$ and $x'_i := y'_i/n$ for $i=1, \ldots, w$.

The basic idea underlying our approach is to find a dissection 
of $B_w$ that gives $T_w$. The encoding and decoding
algorithms are obtained by tracking how the points
$\vec{y}$ and $\vec{y'}$ move during the dissection.

The volume of $B_w$ is $1
\times \frac{1}{2}
\times \frac{1}{3}
\times \cdots
\times \frac{1}{w} = \frac{1}{w!}$.
This is also the volume of $T_w$, as the following argument shows.
Classify the points $\vec{x} = (x_1, \ldots, x_w)$
in the unit cube $[0,1]^w$ into
$w!$ regions according to their order when sorted;
the regions are congruent, so all have volume $1/w!$,
and the region where the $x_i$
are in nondecreasing order is $T_w$.

\begin{figure}[t]
\begin{center}
\includegraphics[width=9cm]{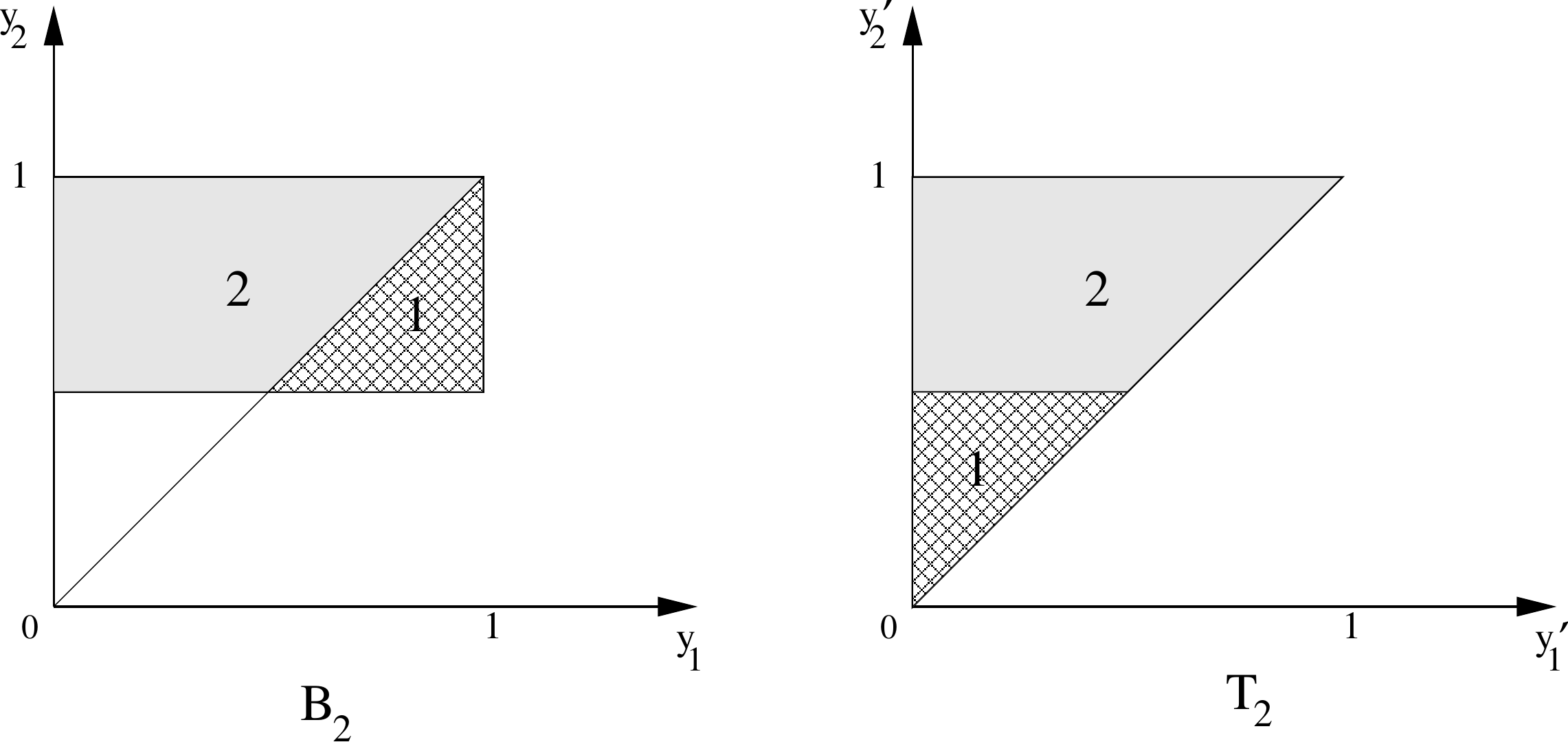}
\end{center}
\caption{ Two ways to dissect rectangle $B_2$ to give triangle $T_2$.
Piece 1 may be rotated about center into its new position,
or reflected in main diagonal and translated downwards. }
\label{fig:2D}
\label{fig2}
\end{figure}

We now return to the case $w=2$. There are many
ways to dissect the rectangle $B_2$ into the right 
triangle $T_2$. We will consider two such dissections,
both two-piece dissections based on Fig.~\ref{fig:2D}.

In the first dissection, the triangular piece
marked 1 in Fig.~\ref{fig:2D} is rotated clockwise
about the center of the square until it reaches the position shown on the right
in Fig.~\ref{fig:2D}. In the second dissection,
the piece marked 1 is first reflected in the main diagonal
of the square and then translated downwards 
until it reaches the position shown on the right
in Fig.~\ref{fig:2D}.
In both dissections the piece marked 2 is fixed.

The two dissections can be specified in
terms of coordinates\footnote{For our use of
a mixture of $\leq$ and $<$ signs, 
see the remark at the end of Section \ref{sec:review}.}
as follows.
For the first dissection, we set
\begin{eqnarray}
\left\{
\begin{array}{ll}
(x'_1,x'_2)  :=  (x_1,x_2) & \mbox{if~} x_1 < x_2 \\
(x'_1,x'_2)  :=  (1-x_1,1-x_2) & \mbox{if~} x_1 \ge x_2 
\end{array}
\right.
\label{eq2Da}
\end{eqnarray}
and for the second, we set
\begin{eqnarray}
\left\{
\begin{array}{ll}
(x'_1,x'_2)  :=  (x_1,x_2) & \mbox{if~} x_1 < x_2 \\
(x'_1,x'_2)  :=  (x_2 - \frac{1}{2},x_1 - \frac{1}{2}) &
  \mbox{if~} x_1 \ge x_2 
\end{array}
\right.
\label{eq2Db}
\end{eqnarray}

The first dissection involves only a rotation, but seems
harder to generalize to higher dimensions.
The second one is the one we will generalize;
it uses a reflection, but as mentioned at
the end of Section \ref{sec:review},
this is permitted by the definition of a dissection.

We next illustrate how these dissections can be converted
into encoding algorithms for constant weight (weight $2$)
binary codes. Again there may be several solutions,
and the best algorithm may depend
on arithmetic properties of $n$ (such as its parity).
We work now with the unnormalized sets
${\cal R}_2$ and ${\cal C}_2$.
In each case the output is a weight-$2$ binary  vector
with $1$'s in positions $y'_1$ and $y'_2$.

\subsection{First Dissection, Algorithm 1}
\begin{enumerate}
\item The input is an information vector $(y_1, y_2) \in {\cal R}_2$ with
$1 \leq  y_1 \leq  n-1$ and $\lceil n/2 \rceil+1 \leq y_2 \leq n$.
\item If $y_1<y_2$, we set $y'_1=y_1$, $y'_2=y_2$, otherwise
we set $y'_1 = n-y_1$ and $y'_2 = n-y_2+1$.
\end{enumerate}
For $n$ even, this algorithm generates all possible $n(n-1)/2$ codewords.
For $n$ odd it generates only $(n-1)^2/2$ codewords,
leading to a slight inefficiency, and the following algorithm
is to be preferred.

\subsection{First Dissection, Algorithm 2}
\begin{enumerate}
\item The input is an information vector $(y_1, y_2) \in {\cal R}_2$ with
$1 \leq y_1 \leq n$, $\lceil (n+1)/2\rceil +1 \leq y_2 \leq n$.
\item If $y_1<y_2$, we set $y'_1=y_1$, $y'_2=y_2$, otherwise
we set $y'_1 = n-y_1+1$, $y'_2 = n-y_2+2$.
\end{enumerate}
For $n$ odd, this algorithm generates all $n(n-1)/2$ codewords,
but for $n$ even it generates only $n(n-1)/2$ codewords,
again leading to a slight inefficiency.

\subsection{Second Dissection}
\begin{enumerate}
\item The input is an information vector $(y_1, y_2) \in {\cal R}_2$ with
$1 \leq  y_1 \leq  n-1$ and $\lceil n/2 \rceil+1 \leq y_2 \leq n$.
\item If $y_1<y_2$, we set $y'_1=y_1$, $y'_2=y_2$, otherwise
we set $y'_1 = y_2 -\lceil n/2 \rceil$,
$y'_2 = y_1 -\lceil n/2 \rceil +1$.
\end{enumerate}
For $n$ even, this algorithm generates all $n(n-1)/2$ codewords,
but for $n$ odd it generates only $(n-1)^2/2$ codewords,
leading to a slight inefficiency. There is a similar
algorithm, not given here, which is better when $n$ is odd.

Note that only one test is required in any of the encoding
algorithms. 
The mappings are  invertible, 
with obvious decoding algorithms corresponding to the
inverse mappings from 
${\cal C}_2$ to ${\cal R}_2$

We now extend this method to weight $w=3$. 
Fortunately, the Dehn invariants for both the brick $B_3$ and our 
unit orthoscheme $T_3$,
which is the tetrahedron\footnote{To solve Hilbert's third problem,
Dehn showed that this tetrahedron is not equidecomposable
with a regular tetrahedron of the same volume.}
with vertices $(0,0,0), (0,0,1), (0,1,1)$ and $(1,1,1)$,
are zero (since in both cases all dihedral angles are 
rational multiples of $\pi$), and so by the Dehn-Sydler theorem
the polyhedra $B_3$ and $T_3$ {\em are} equidecomposable.
As already mentioned in Section \ref{sec:intro}, the Dehn-Sydler theorem
applies only in three dimensions. But it will follow from the algorithm
given in the next section
that $B_w$ and $T_w$ are equidecomposable in all dimensions.

\begin{figure}[t]
\begin{center}
\includegraphics[width=7cm]{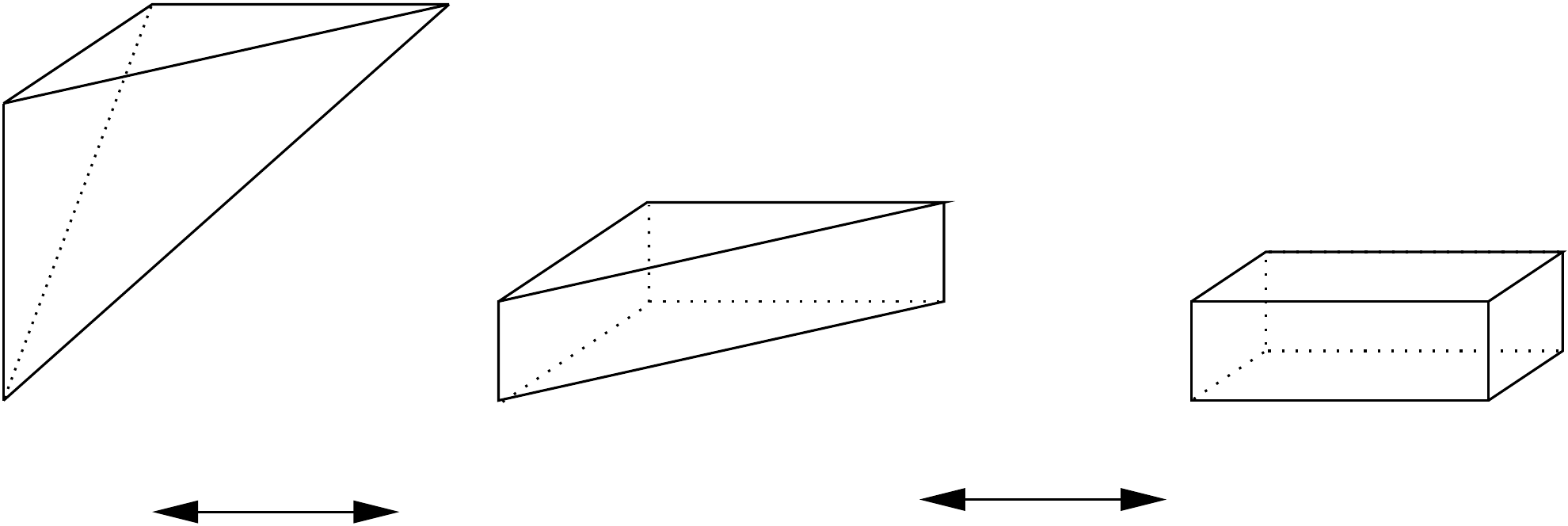}
\end{center}
\caption{Transformation from tetrahedron to rectangular prism.}
\label{fig:3D}
\label{fig3}
\end{figure}

We will continue to describe
the encoding step 
(the map from $B_w$ to $T_w$) first.
We will give an inductive dissection (see Fig. \ref{fig:3D}),
transforming $B_3$ to $T_3$ in two steps, effectively
reducing the dimension by one at each step. 
In the first step, the brick $B_3$ is dissected
into a  triangular prism (the product of a right triangle, $T_2$, and an interval),
and in the second step this triangular prism is dissected into the tetrahedron $T_3$.
Note that the first step has essentially been solved 
by the dissection given in Eqn. (\ref{eq2Db}).

\begin{figure}[t]
\begin{center}
\includegraphics[width=6cm]{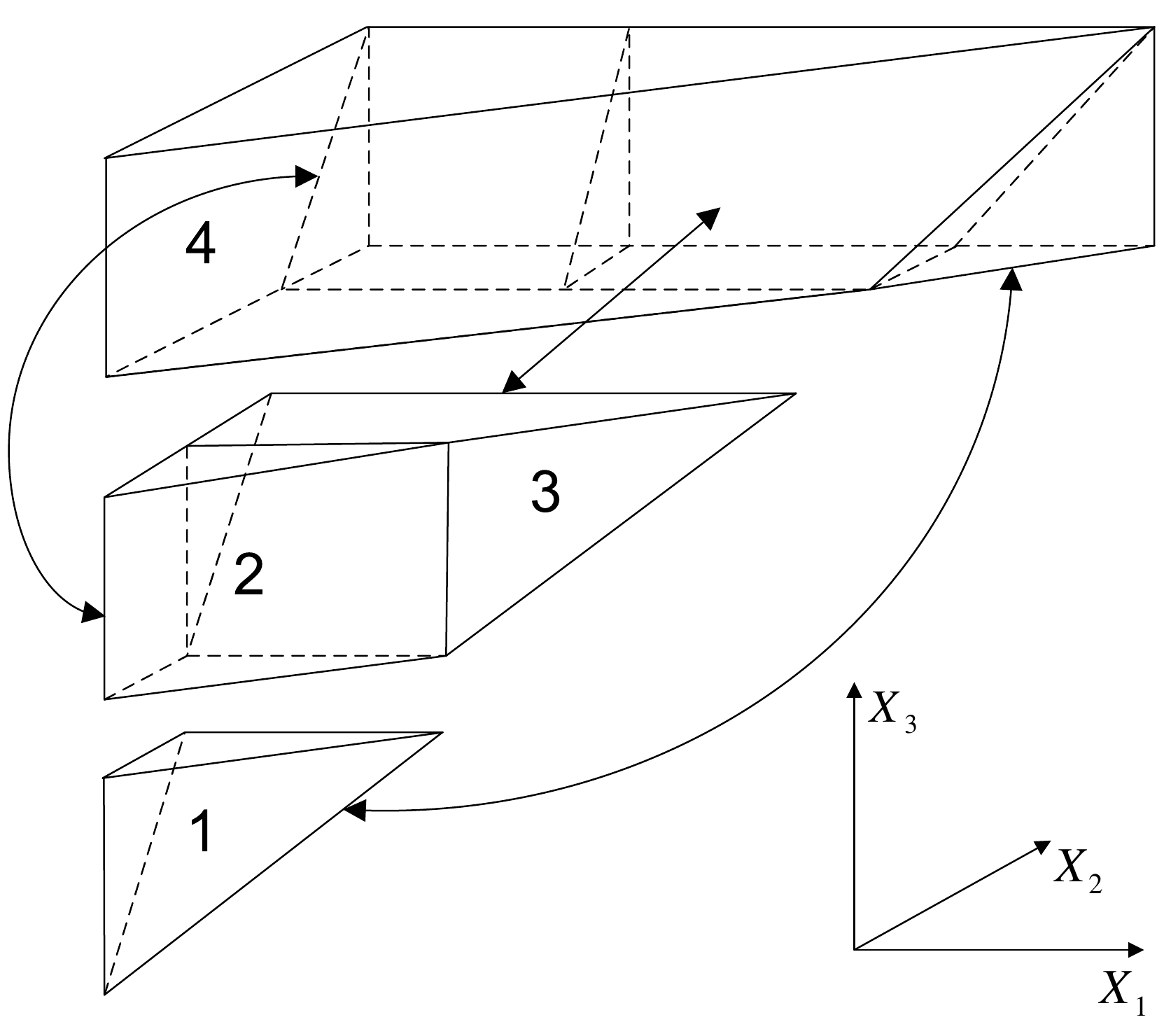}
\end{center}
\caption{Four-piece dissection of tetrahedron to triangular prism.
Pieces 2 and 3 are reflected.}
\label{fig:3Dpiece}
\label{fig4}
\end{figure}

For the second step we use a four-piece dissection of 
the triangular prism to the tetrahedron $T_3$.
This dissection, shown with the tetrahedron and prism superimposed 
in Fig. \ref{fig:3Dpiece}, appears to be new.

There is a well-known dissection of the same pair of polyhedra
that was first published by Hill in 1896~\cite{Hill1896}.
This also uses four pieces, and is discussed in several references: see
Boltianskii~\cite[p.~99]{Bolt1978},
Cromwell~\cite[p.~47]{Crom1997},
Frederickson~\cite[Fig.~20.4]{Fred1997},
Sydler~\cite{Sydl1956},
Wells~\cite[p.~251]{Well1991}.
However, Hill's dissection seems harder to generalize to higher dimensions.
Hill's dissection does have the advantage over ours
that it can be accomplished purely by translations and rotations,
whereas in our dissection two of the pieces (pieces labeled
2 and 3 in Fig. \ref{fig:3Dpiece}) 
are also reflected. However, as mentioned at
the end of Section \ref{sec:review},
this is permitted by the definition of a dissection,
and is not a drawback for our application. \footnote{This
dissection would also work if piece 2 was merely translated and rotated,
not reflected, but the reflection is required by our
general algorithm.}
Apart from this, our dissection is simpler than Hill's,
in the sense that his dissection requires a cut along
a skew plane ($x_1-x_3=1/3$),
whereas all our cuts are parallel to coordinate axes.

To obtain the four pieces shown
in Fig. \ref{fig:3Dpiece}, we first make
two horizontal cuts along the planes
$x_3 = \frac{1}{3}$ and $x_3 = \frac{2}{3}$, dividing the 
tetrahedron into three slices. We then cut the middle slice into two 
by a vertical cut along the plane $x_2 = \frac{1}{2}$.

There appears to be a tradition in geometry books 
that discuss dissections  of not giving coordinates for the pieces.
To an engineer this seems unsatisfactory,
and so in Table \ref{T1} we list the vertices of
the four pieces in our dissection.
Piece 1 has four vertices, while the other three
pieces each have six vertices.
(In the Hill dissection the numbers of vertices
of the four pieces are $4$, $5$, $6$ and $6$ respectively.)
Given these coordinates,
it is not difficult to verify that the four
pieces can be reassembled to form the triangular prism,
as indicated in Fig. \ref{fig:3Dpiece}. As
already remarked, pieces 2 and 3 
are also reflected (or ``turned over'' in a fourth dimension).
The correctness of the dissection
also follows from the alternative description of this
dissection given below.

\renewcommand{\arraystretch}{1.3}

\begin{table}[htb]
$$
\begin{array}{|c|l|} \hline
\mbox{Piece} & \mbox{Coordinates} \\ \hline
1 & [0,0, 0], [0,0,1/3], [0,1/3,1/3], [1/3, 1/3,1/3]. \\
2 & [0,0,1/3], [0,1/3,1/3], [1/3,1/3,1/3], \\
  & [0,0,2/3], [0,1/3,2/3], [1/3,1/3,2/3]. \\
3 & [0,1/3,1/3], [1/3,1/3,1/3], [0, 1/3,2/3], \\
  & [0,2/3,2/3], [2/3,2/3,2/3], [1/3,1/3,2/3]. \\
4 & [0,0,2/3], [0,2/3,2/3], [2/3, 2/3,2/3], \\
  & [0,0,1], [0,1,1], [1, 1,1]. \\
\hline
\end{array}
$$
\caption{ Coordinates of vertices of pieces in dissection
of tetrahedron shown in Fig. \ref{fig:3Dpiece}.  }
\label{T1}
\end{table}

\renewcommand{\arraystretch}{1.0}

The dissection shown in Fig. \ref{fig:3Dpiece} can be described
algebraically as follows. We describe it in the more
logical direction, going from the triangular prism to
the tetrahedron 
since this is what we will generalize to
higher dimensions in the next section.
The input is a vector $(x_1, x_2, x_3)$ with
$0 \le x_1 \le x_2 <1$, $\frac{2}{3} \le x_3 <1$;
the output is a vector $(x'_1, x'_2, x'_3)$ with
$0 \le x'_1 \le x'_2 \le x'_3 <1$, given by $(x'_1, x'_2, x'_3) =$
\renewcommand{\arraystretch}{1.3}
\begin{eqnarray}
\left\{
\begin{array}{ll}
(x_1,x_2,x_3)       & \mbox{if~} x_1 \le x_2 <x_3 \\
(x_1-\frac{1}{3},x_3-\frac{1}{3},x_2-\frac{1}{3})  &
\mbox{if~}\frac{1}{3} \le x_1 <x_3 \le x_2\\
(x_3-\frac{2}{3},x_2-\frac{2}{3},x_1+\frac{1}{3})  &
\mbox{if~}x_1 \le \frac{1}{3} <x_3 \le x_2\\
(x_3-\frac{2}{3},x_1-\frac{2}{3},x_2-\frac{2}{3})  &
\mbox{if~}x_3 \le x_1 \le x_2
\end{array}
\right.
\label{Eqf3}
\end{eqnarray}
\renewcommand{\arraystretch}{1.0}

The four cases in Eqn. (\ref{Eqf3}), after being transformed, correspond to the pieces labeled
4, 3, 2, 1 respectively in Fig. \ref{fig:3Dpiece}. We see from  Eqn. (\ref{Eqf3})
that in the second and third cases the linear
transformation has determinant $-1$, indicating
that these two pieces must be reflected.

\begin{figure*}[th]
\begin{center}
\includegraphics[width=14cm]{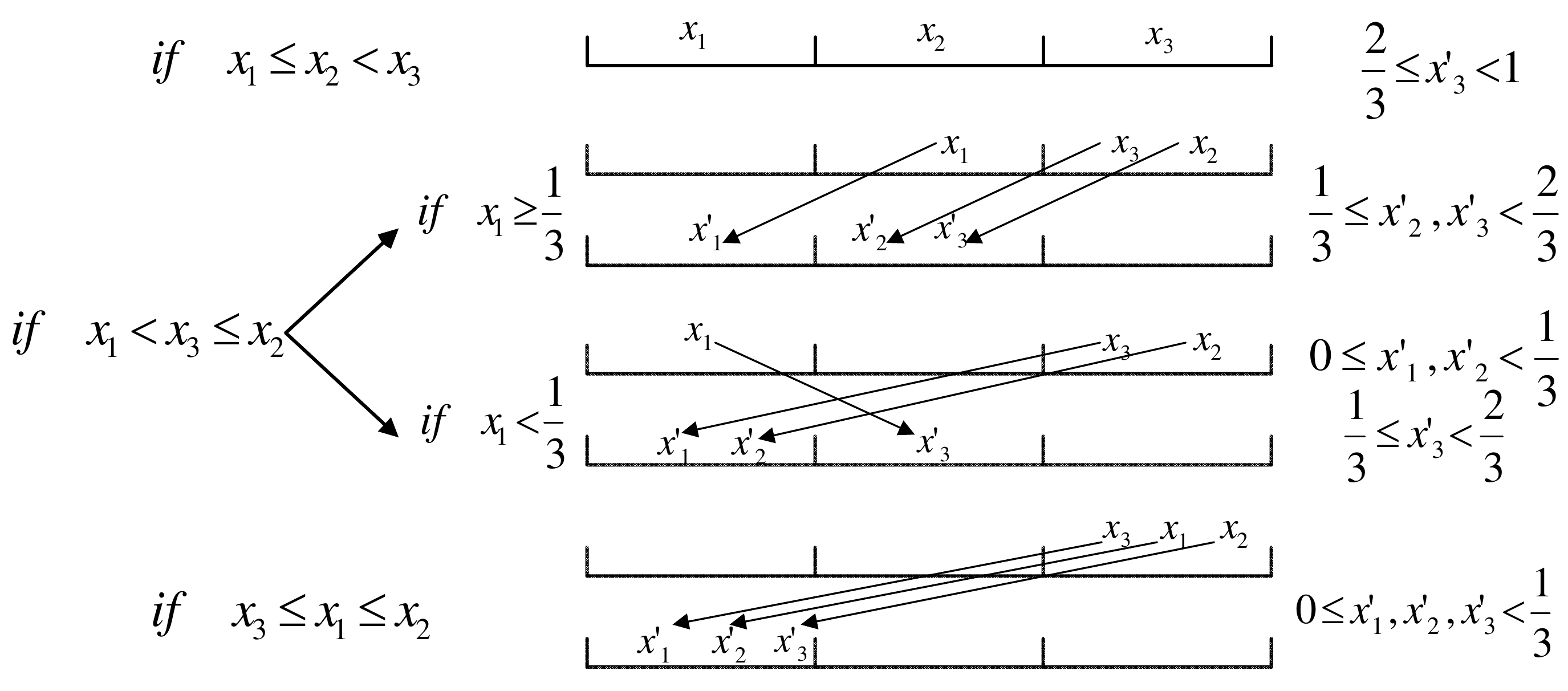}
\end{center}
\caption{Cut-and-paste description of the inverse transformation
from triangular prism to tetrahedron.}
\label{fig:shift}
\label{fig5}
\end{figure*}

Since it is hard to visualize dissections in dimensions greater than three,
we give a schematic representation of
the above dissection that avoids drawing polyhedra.
Fig.~\ref{fig:shift} shows a representation of the transformation 
from the triangular prism  to
the tetrahedron  $T_3$,
equivalent to that given in  Eqn. (\ref{Eqf3}).
The steps shown in Fig.~\ref{fig:shift} may be referred to as 
``cut and paste'' operations,
because, as Fig.~\ref{fig:shift} shows,
the vector in the triangular prism is literally cut up into pieces
which are rearranged and relabeled.
Note that, to complete the transformation,
we precede this operation by
the dissection given in Eqn. (\ref{eq2Db}),
finally establishing the bijection between $B_3$ and $T_3$.

We now describe the mapping shown in Fig.~\ref{fig:shift}
in more detail. The triangular prism is represented by the
set of partially ordered triples $(x_1,x_2,x_3)$ with
$0\leq x_1 \leq x_2<1$ 
and $\frac{2}{3}\leq x_3<1$, and we wish to transform
this into the tetrahedron consisting of the points
$(x'_1, x'_2, x'_3)$ with $0\leq x'_1 \leq x'_2 \leq x'_3 <1$.

We divide the interval $[0,1)$ into $w=3$ equal segments
of length $1/w = 1/3$, and consider where the
points $x_1, x_2$ and $x_3$ fall in this interval,
given that $(x_1, x_2, x_3)$ is in the triangular prism.
There are three possibilities for where $x_3$ lies in relation to 
$0 \le x_1 \le x_2 < 1$, and we further divide
the case $x_1 \le x_3 <x_2$ into two
subcases depending on whether $x_1 \ge \frac{1}{3}$ or $x_1 < \frac{1}{3}$.
These are the four cases shown in Fig.~\ref{fig:shift},
and correspond one-to-one with the four
dissection pieces in Fig.~\ref{fig:3Dpiece}. 
Fig.~\ref{fig:shift} shows how the triples 
$x_1, x_2, x_3$  (reindexed according to their
relative positions) are mapped to the triples 
$x'_1, x'_2, x'_3$.

The last column of Fig.~\ref{fig:shift} shows the ranges of 
the $x'_i$ in the four cases; the fact that these
ranges are disjoint guarantees that the mapping
from $x_1, x_2, x_3$ to $x'_1, x'_2, x'_3$ is invertible.
The ranges of the $x'_i$ will be discussed in more detail in the
following section after the general algorithms are presented. 

This operation can now be described without explicitly 
mentioning the underlying  dissection.  Each interval of
length $1/w$, together with the given $x_i$ values 
within it, is treated as a single complete unit.
In the  ``cut and paste" operations, these units are rearranged
and relabeled in such a way that the operation is invertible.

\section{Algorithms and Proof of Correctness}
\label{sec:alg}

In the previous section we provided an encoding and decoding algorithm
for weights $w=2$ and $w=3$,
based on our geometric interpretation of ${\cal C}_2$ and ${\cal C}_3$
as points in $\Reals^w$.
In this section, the algorithm is 
generalized to larger values of the weight $w$. 
We start with the geometry, 
and  give a dissection 
of the ``brick'' $B_w$ into the orthoscheme $T_w$.
We work with the normalized coordinates
$x_i=y_i/n$ (for a point in $B_w$) and $x'_i=y'_i/n$ 
(for a point in $T_w$), where $1\leq i \leq w$. 
Later in this section,
we discuss the modifications needed to take
into account the fact that the $y'_i$ must be integers.

\subsection{An Inductive Decomposition of the Orthoscheme}

Restating the problem, we wish
to find a bijection  $F_w$ between the sets $B_w$ and $T_w$. 
The inductive approach developed for $w=3$ 
(where the $w=2$ case was a subproblem) will be generalized. 
Of course the bijection $F_1$ between $B_1$ and $T_1$ is trivial.
We assume that a bijection $F_{w-1}$ is known
between $B_{w-1}$ and  $T_{w-1}$, 
and show how to construct a bijection $F_w$ between $B_w$ and $T_w$. 

The last step in the induction uses a map $f_w$
from the prism $T_{w-1}\times[1-\frac{1}{w},1)$ to ${T_w}$
($f_2$ is the map described in Eqn. (\ref{eq2Db})
and $f_3$ is described in Eqn. (\ref{Eqf3})).
The mapping $F_w$ from $B_w$ to $T_w$ is then given recursively by
$F_w: (x_1,x_2,\ldots,x_w) \mapsto (x'_1,x'_2,\ldots,x'_w)$, where
\begin{equation}\label{EQf}
(x'_1,x'_2,\ldots,x'_w) := f_w(F_{w-1}(x_1,x_2,\ldots,x_{w-1}),x_w) \,.
\end{equation}
For $w=1$ we set
$$
F_1  :=  f_1: B_1 \rightarrow T_1, ~ (x_1) \mapsto (x'_1) = (x_1) \,.
$$ 
By iterating Eqn. (\ref{EQf}), we see that $F_w$ is obtained by
successively applying the maps $f_1, f_2, $\ldots, $f_w$.

The following algorithm defines $f_w$ for $w \ge 2$.
We begin with an algebraic definition of
the mapping and its inverse, and
then discuss it further in the following section.
The input to the mapping $f_w$ is a vector $\vec{x} := (x_1,x_2,\ldots,x_w)$,
with $(x_1,x_2,\ldots,x_{w-1}) \in T_{w-1}$ and $x_w \in [1-1/w,1)$;
the output is a vector $\vec{x'} := (x'_1,x'_2,\ldots,x'_w) \in T_w$. 

\begin{flushleft}
{\em Forward mapping $f_w$} ($w \ge 2$):
\end{flushleft}

1) Let 
\begin{eqnarray}
i_0& := &\min \{ i \in \{1,\ldots,w\} \mid x_w \le x_i \} \,. \nonumber
\end{eqnarray}

2) Let
\begin{eqnarray}
j_0& := &\min \{ i  \in \{1,\ldots, i_0\} \mid
  w-i_0+i-1 \le w x_i \} -1 \,. \nonumber
\end{eqnarray}

3) Set $x_k^{\prime}$ equal to:
\begin{eqnarray}
\left\{
\begin{array}{ll}
x_{k+j_0}-\frac{w+j_0-i_0}{w}       & \mbox{for~}k=1,\ldots,i_0-j_0-1\\
x_w-\frac{w+j_0-i_0}{w}             & \mbox{for~}k=i_0-j_0        \\
x_{k+j_0-1}-\frac{w+j_0-i_0}{w}     & \mbox{for~}k=i_0-j_0+1,\ldots,w-j_0\\
x_{k-w+j_0}+\frac{i_0-j_0}{w}     & \mbox{for~}k=w-j_0+1,\ldots,w   
\end{array}
\right. 
\label{eqn-f-shift}
\end{eqnarray}
Eqn. (\ref{eqn-f-shift}) identifies the ``cut and paste" operations required to
obtain $x'_k$ for different ranges  of the variable $k$.
If the initial index in one of the four cases in Eqn. (\ref{eqn-f-shift})
is smaller than the final index, that case is to be skipped.
A case is also skipped if the subscript
for an $x_i$ is not in the range $1,\ldots,w$.
Note in Step 1 that $i_0 = w$ if $x_w$ is the
largest of the $x_i$'s.
This implies that $j_0=0$, and then Step 3 is the identity map.

The inverse mapping $G_w$ from $T_w$ to $B_w$
has a similar recursive definition.
The $w$th step in the induction is the map 
$g_w: T_w \rightarrow {T_{w-1}\times[1-\frac{1}{w},1)}$
defined below.
For $w=1$ we set
$$
G_1  :=  g_1: T_1 \rightarrow B_1, (x'_1) \mapsto (x_1) = (x'_1) \,.
$$
The map $G_w$ is obtained by
successively applying the maps $g_w, g_{w-1}, $\ldots, $g_1$.

\begin{flushleft}
{\em Inverse mapping $g_w$} ($w \ge 2$):
\end{flushleft}

1) Let 
\begin{eqnarray}
m_0& := &\max \{ i \in \{1,\ldots,w\} \mid i-1 \le w x'_i \} \,. \nonumber
\end{eqnarray}

2) If $m_0=w$, let $j_0 := 0$,
otherwise let
\begin{eqnarray}
j_0& := &w- \max \{ i  \in \{m_0 + 1,\ldots, w\} \mid
  wx'_i \le m_0 \} \,; \nonumber
\end{eqnarray}
in either case, let $i_0 := j_0 + m_0$.

3) Set $x_k$ equal to:
\begin{eqnarray}
\left\{
\begin{array}{ll}
x_{k+w-j_0}^\prime-\frac{i_0-j_0}{w}  & \mbox{for~} k=1,\ldots,j_0 \\
x_{k-j_0}^\prime+\frac{w+j_0-i_0}{w}        & \mbox{for~} k=j_0+1,\ldots,i_0-1\\
x_{k-j_0+1}^\prime+\frac{w+j_0-i_0}{w}      & \mbox{for~} k=i_0,\ldots,w-1\\
x_{i_0-j_0}^\prime+\frac{w+j_0-i_0}{w}      & \mbox{for~} k=w
\end{array}
\right.
\label{eqn-g-shift}
\end{eqnarray}

Note that the transformations in Eqn. (\ref{eqn-f-shift})
and Eqn. (\ref{eqn-g-shift}) are formal inverses of each other,
and that these transformations are volume-preserving.
The underlying linear transformations are orthogonal 
transformations with determinant $+1$ or $-1$.

Before proceeding further, let us verify that in the case $w=3$,   
the mapping $f_w = f_3$ agrees with that given in Eqn. (\ref{Eqf3}).
\begin{itemize}
\item[$\bullet$]If $x_1\leq x_2<x_3$, then $i_0=3$, $j_0=0$ and the map is the
identity, as mentioned above.
\item[$\bullet$]If $x_1<x_3\leq x_2$ there are two subcases:
\begin{itemize}\item[$\circ$]If $\frac{1}{3} \le x_1$ then $i_0=2$, $j_0=0$.
\item[$\circ$]If $x_1<\frac{1}{3}$ then $i_0=2$, $j_0=1$.
\end{itemize}
\item[$\bullet$]If $x_3\leq x_1\leq x_2$, then $i_0=1$, $j_0=0$.
\end{itemize}
The transformations in Eqn. (\ref{eqn-f-shift}) now exactly
match those in Eqn. (\ref{Eqf3}).

\subsection{Interpretations and Explanations}

\begin{figure*}[htb]
\begin{center}
\includegraphics[width=14cm]{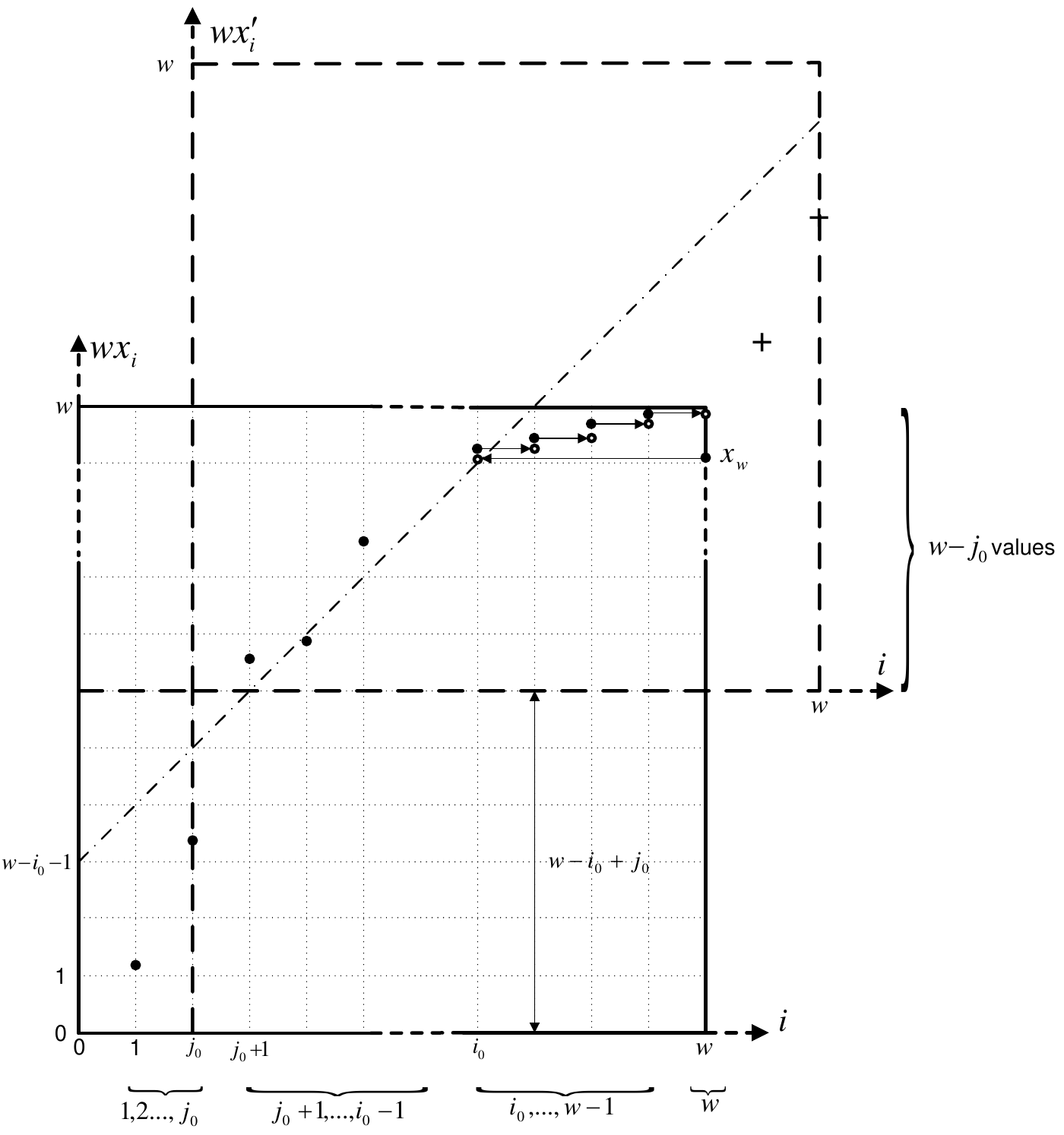}
\end{center}
\caption{A graphical illustration of the forward and inverse mapping.}
\label{fig:intuitive}
\end{figure*}

In Fig. \ref{fig:intuitive}, we give a graphical interpretation of the
algorithm, which can be regarded as a 
generalization of the ``cut and paste'' description given above. 
This figure shows the transformation defined by the $w$th step $f_w$
in the algorithm.
At this step, we begin with a list of 
$w-1$ numbers $(x_1,x_2,\ldots,x_{w-1})$  in
increasing order, and a further number $x_w$ which 
may be anywhere in the interval $[1-1/w,1)$. 
This list of $w$ numbers is plotted
in the plane as the set of $w$ points $(i,wx_i)$ for $i=1,2,\ldots,w$
(indicated by the solid black circles in  Fig. \ref{fig:intuitive}). 
In the first step in the forward algorithm, 
the augmented list $(x_1,x_2,\ldots,x_w)$
is sorted into increasing order.  
In the sorted list, $x_w$ now occupies position ${i_0}$,
so the point $(w, w x_w)$ moves to the left,
to the new position $(i_0,wx_w)$, and the points
$(i, wx_i)$ for $i=i_0+1,\ldots,w-1$ move to the right.
This is indicated by the arrows in the figure. 
The new positions of these points are
marked by hollow circles.

The point $(i_0,wx_w)$ now lies 
between the grid points $(i_0,w)$ and $(i_0,w-1)$
(it may coincide with the latter point),
since $x_w \ge 1-\frac{1}{w}$.
We draw the line $y=x+w-i_0-1$
(shown as the dashed-and-dotted line in Fig. \ref{fig:intuitive}).
This has unit slope and passes through the points
$(i_0,w-1)$ and $(0,w-i_0-1)$.
The algorithm then computes $j_0+1$ to be the smallest index $i$
for which $x_{i}$ is on or above this line.
Once $i_0$ and $j_0$ have been determined, 
the forward mapping proceeds as follows.
The points $(i, wx_i)$ for $i=1,\ldots,j_0$
are shifted to the right of the figure
and are moved upwards by the amount $(i_0-j_0)/w$,
their new positions being indicated by crosses in
the figure.
Finally, the origin is moved to the grid point
$(j_0,w-i_0+j_0)$ and the
points are reindexed.
The $m_0 := i_0-j_0$ points
which originally had indices $j_0+1, \ldots,i_0$
become points $1,\ldots,m_0$ after reindexing.
In the new coordinates, the final positions
of the points lie inside the square
region $[1,w)\times[1,w)$.
The reader can check that this process is
exactly equivalent to the algebraic description of $f_w$ given above.

To recover $i_0$ and $j_0$, we first determine the value of $m_0 := i_0-j_0$. 
This can indeed be done since $m_0$ is precisely the index
of the largest $wx'_i$ that lies on or above the  line $y=x-1$
in the new coordinate system. 
Note that the position of this  line is 
independent of $i_0$ and $j_0$ and $(x'_1,x'_2,\ldots,x'_w)$.
This works because the points $wx_1, \ldots, wx_{j_0}$ in
the original coordinate system, before the origin
is shifted, are moved right by $w$ units
and upwards by $w$ units,
so points below the dashed-and-dotted line
remain below the line.
Furthermore, observe that in the new coordinate system
the number of points $(i,wx'_i)$ below the line $y=m_0$
is equal to $w-j_0$. 
Thus the correct $i_0$ and $j_0$ values may be recovered,
and the inverse mapping can be successfully performed.

The following remarks record two properties of the algorithm that will
be used later.

{\em Remark 1:}
Step 2 of the forward algorithm implies that
$x_{j_0}<\frac{w-i_0+j_0-1}{w}$ and
$x_{j_0+1}\geq{\frac{w-i_0+j_0}{w}}$.
It follows that there is no $i$ in the range
$1\leq{i}\leq{w}$ for which
$$
\frac{w-i_0+j_0-1}{w}\leq{x_i}<\frac{w-i_0+j_0}{w} \,.
$$

{\em Remark 2:} The forward algorithm produces a vector $\vec{x'}$ whose components satisfy
\begin{eqnarray}
\label{eqn:moveleft}
0\leq{x_1^\prime}\leq \cdots
\leq{x_{i_0-j_0}^\prime}\leq \cdots \leq{x_{w-j_0}^\prime}<\frac{i_0-j_0}{w} \,,
\end{eqnarray}
\begin{eqnarray}
\label{eqn:moveright}
\frac{i_0-j_0}{w}\leq{x_{w-j_0+1}^\prime} \le
{x_{w-j_0+2}^\prime} \le \cdots \le {x_{w}^\prime} <1 \,,
\end{eqnarray}
and 
\begin{eqnarray}
\label{Eqbelow}
{x_{k}^\prime}<\frac{k-1}{w}, \mbox{~for~}  w-j_0+1\leq{k}\leq{w} \,.
\end{eqnarray}
Eqns. (\ref{eqn:moveleft}) and (\ref{eqn:moveright}) follow from the
minimizations in Steps~1 and 2 of the forward algorithm, respectively. 
The right-hand side of Eqn. (\ref{Eqbelow}) expresses the fact, already mentioned,
that the first $j_0$ points remain below the dotted-and-dashed line
after they are shifted.

\subsection{Proof of Correctness}

We now give the formal proof that the algorithm is correct.
This is simply a matter of collecting together
facts that we have already observed. 

\begin{theorem}\label{Th1}
For any $w \ge 1$, the forward mapping $f_w$  is a one-to-one mapping
from ${T_{w-1}\times{[1-\frac{1}{w},1)}}$ to $T_w$ with inverse $g_w$.
\end{theorem}

\begin{myproof}
First, it follows from Remark 2 that, for 
$\vec{x} \in  {T_{w-1}\times{[1-\frac{1}{w},1)}}$,
$\vec{x'}=(x_1', x_2', \ldots,x_w')$ satisfies
$0\leq{x_1'}\leq{x_2'}\leq \cdots \leq{x_w'}<1$, and so
is an element of $T_w$.

Suppose there were two different choices for $\vec{x}$, say
$\vec{x}^{(1)}$ and $\vec{x}^{(2)}$, such that
$$
f_w(\vec{x}^{(1)}) = f_w(\vec{x}^{(2)}) = \vec{x}' \,.
$$
We know that $\vec{x}'$ determines $m_0, j_0$ and $i_0$.
So $\vec{x}^{(1)}$ and $\vec{x}^{(2)}$ have the same associated values
of $i_0$ and $j_0$. 
But for a given pair $(i_0,j_0)$, Eqn. (\ref{eqn-f-shift}) is invertible.
Hence  $\vec{x}^{(1)} = \vec{x}^{(2)}$, and $f_w$ is one-to-one.

Note that the transformations in Eqn. (\ref{eqn-f-shift})
and Eqn. (\ref{eqn-g-shift}) are inverses of each other. Hence
$f_w$ is also an onto map, and $g_w$ is its inverse.
\end{myproof}

\subsection{Number of Pieces}
The map $f_w$, which dissects
the prism $T_{w-1}\times[1-\frac{1}{w},1)$ to give the
orthoscheme ${T_w}$, has one piece for each pair $(i_0, j_0)$.
If $i_0=w$ then $j_0=0$,
while if $1 \le i_0 \le w-1$, $j_0$ takes all values
from $0$ to $i_0-1$.
(It is easy to write down an explicit point in the interior of 
the piece corresponding to a specified pair of values
of $i_0$ and $j_0$. Assume $i_0 <w$ and
set $\delta = 1/w^3$.
Take the point with coordinates
$(x_1, \ldots,x_w)$ given by
$x_w = (w-1)/w+\delta$;
$x_i = x_w + \delta(i-i_0)$ 
for $i=i_0+1, \ldots, w-1$;
$x_i = (i+w-i_0-1-\delta)/w$
for $i=1, \ldots,j_0$;
$x_i = (i+w-i_0-1+\delta)/w$
for $i=j_0+1, \ldots,i_0-1$.)
The total number of
pieces in the dissection is therefore
$$
1+1+2+3+\cdots+(w-1) = \frac{w^2-w+2}{2} \,,
$$
which is $1,2,4,7,11,\ldots$ for $w=1,2,3,4,5,\ldots$.
This is a well-known sequence, entry A124 in \cite{OEIS},
which by coincidence also arises in a different
dissection problem: it is the maximal number
of pieces into which a circular disk can be cut
with $w-1$ straight cuts.
For example, with three cuts, a pizza can be
cut into a maximum of seven pieces, and
this is also the number of pieces in the
dissection defined by $f_4$.

\subsection{The Algorithms for Positive Integers}
To apply the above algorithm to the problem of
encoding and decoding constant weight codes,
we must work with 
positive integers rather than real numbers,  which entails a
certain loss in rate, although
the algorithms remain largely unchanged. 
Let $\NN := \{1,2,3,\ldots\}$, and let $n$ and $w$ be given with $2w<n$.
In a manner analogous to the real-valued case,
we  find a bijection between 
a finite hyper-rectangle or brick
$B_w^{\NN}\subset{{\NN}^w}$ and a subset of the finite orthoscheme
$T_w^{\NN}\subset{{\NN}^w}$,
where 
$B_w^{\NN}$ is the set of vectors $(y_1,y_2,\ldots,y_w)\in{{\NN}^w}$ satisfying
$$
n-(w-i)-\lfloor{\frac{n-(w-i)}{i}}\rfloor+1\leq{y_i}\leq{n-(w-i)} \,, 
$$
for $i=1,2,\ldots,w$,
and 
$T_w^{\NN}$ is the set of vectors $(y_1,y_2,\ldots,y_w)\in{{\NN}^w}$ satisfying
$$
1\leq{y_1}<y_2< \cdots <y_w\leq{n} \,.
$$ 
Note that usually $|B_w^{\NN}| < |T_w^{\NN}|$, which entails a loss in rate.

The forward mapping $f_w$ is now replaced by the map
$f^{\NN}_w$,
which sends 
$(y_1,y_2,\ldots,y_w)$ with $(y_1,y_2,\ldots,y_{w-1})\in{T_{w-1}^{\NN}}$ 
and
${n-\lfloor\frac{n}{w}\rfloor+1}\leq{y_w}\leq{n}$
to an element of ${T_{w}^{\NN}}$. 
Let us write
$n=pw+q$, where $p\geq{0}$ and $0\leq{q}\leq{w-1}$. 
We partition the range $1,2,\ldots,n$ into $w$ parts,
where the first $n-w-1$ parts each have $p$ elements, 
the next $q$ parts each have $p+1$ elements,
and the last part has $p$ elements (giving a total of $n$ elements).
This is similar to the real-valued case, where each interval had length $1/w$. 

1) Let 
\begin{eqnarray}
i_0& := &\min \{ i \in \{1,\ldots,w\} \mid y_w \le y_i \} \,. \nonumber
\end{eqnarray}

2) Let
\begin{eqnarray}
j_0& := &\min \{ i  \in \{1,\ldots, i_0\} \mid
  V_i < y_i \} -1 \,, \nonumber
\end{eqnarray}
where $V_i := (w-i_0+i-1)p+\max\{q-i_0+i,0\}$.

3) Set $y_k^{\prime}$ equal to:
\begin{eqnarray}
\left\{
\begin{array}{ll}
y_{k+j_0}-V_{j_0+1}        & \mbox{for~} k=1,\ldots,i_0-j_0-1\\
y_w-V_{j_0+1}              & \mbox{for~} k=i_0-j_0        \\
y_{k+j_0-1}+1-V_{j_0+1}      & \mbox{for~} k=i_0-j_0+1,\ldots,w-j_0\\
y_{k-w+j_0}+n-V_{j_0+1}  & \mbox{for~} k=w-j_0+1,\ldots,w   
\end{array}
\right.
\label{eqn-fN-shift}
\end{eqnarray}

The inverse mapping $g_w$ is similarly replaced by the map
$g^{\NN}_w:{T_{w}^{\NN}}\rightarrow\{(y_1,y_2,\ldots,y_w):(y_1,y_2,\ldots,y_{w-1})\in{T_{w-1}^{\NN}},$
 ${n-\lfloor\frac{n}{w}\rfloor+1}\leq{y_w}\leq{n}\}$,
defined as follows. Again, assume $n=pw+q$.

1) Let 
\begin{eqnarray}
m_0& := &\max \{ i \in \{1,\ldots,w\} \mid W_i <  y'_i \} \,, \nonumber
\end{eqnarray}
where $W_i := q+(i-1)p+\min\{i-q-1,0\}$.

2) If $m_0=w$, let $j_0 := 0$,
otherwise let
\begin{eqnarray}
j_0& :=  & w - \max \{ i  \in \{m_0 + 1,\ldots, w\} \mid
  y'_i \le W_{m_0}+p \} \,; \nonumber
\end{eqnarray}
in either case, let $i_0 := j_0 + m_0$.

3) Set $y_k$ equal to:
\begin{eqnarray}
\left\{
\begin{array}{ll}
y_{k+w-j_0}^\prime-p-W_{m_0}      & \mbox{for~} k=1,\ldots,j_0 \\
y_{k-j_0}^\prime+n-p-W_{m_0}        & \mbox{for~} k=j_0+1,\ldots,i_0-1\\  
y_{k-j_0+1}^\prime-1+n-p-W_{m_0}    & \mbox{for~} k=i_0,\ldots,w-1\\
y_{i_0-j_0}^\prime+n-p-W_{m_0}      & \mbox{for~} k=w        
\end{array}
\right.
\label{eqn-gN-shift}
\end{eqnarray}

We omit the proofs, since they are similar to those for the real-valued case.

\subsection{Comments on the Algorithm}

The overall complexity of the transform algorithm is $O(w^2)$, because at each
induction step the complexity is linear in the weight at that step. Recall that the
complexities of the arithmetic coding method and Knuth's complementation method are
both $O(n)$. Thus when the weight $w$ is larger than  $\sqrt{n}$, the geometric
approach is less competitive. When the weight is low, the 
proposed geometric
technique is more efficient, because Knuth's complementation method is not
applicable, while the dissection operations of the proposed algorithm makes it
faster than the arithmetic coding method. Furthermore, due to the structure of the
algorithm, it is possible to parallelize part of the computation within each
induction step to further reduce the computation time. 

So far little has been said about mapping a binary sequence to an integer sequence
$y_1,y_2,\ldots,y_w$ such that $y_i\in[L_i,U_i]$, where $L_i$ and $U_i$ are the lower
and upper bound of the valid range as specified by the algorithm. A straightforward
method is to treat the binary sequence as an integer number and then use ``quotient
and remainder" method to find such a mapping. However, this requires a division
operation, and when the binary sequence is long, the computation is not very
efficient. A simplification is to partition the binary sequence into short
sequences, and map each short binary sequence to a pair of integers, as in the case
of a weight two constant weight codes. Through proper pairing of the ranges, the
loss in the rate can be minimized. 

The overall rate loss has two components,
the first from the rounding involved in using natural numbers,
the second from the loss in the above simplified translation step.
However, when the weight is on the order of $\sqrt{n}$, and $n$ is in the range of
$100-1000$, the rate loss is usually $1-3$ bits per block. For example, when
$n=529$, $w=23$, then the rate loss is 2 bits/block 
compared to the best possible code which would encode $k_0=132$ information bits. 

\section{Conclusion}
\label{sec:con}

We propose a novel algorithm for encoding and decoding constant weight binary codes,
based on dissecting the polytope defined by the
set of all binary words of length $n$ and weight $w$,
and reassembling the pieces to form a hyper-rectangle 
corresponding to the input data.
The algorithm has a natural recursive structure, which enables us to 
give an inductive proof of its correctness. 
The proposed algorithm has complexity $O(w^2)$,
independent of the length of the codewords $n$. 
It is especially suitable for constant weight codes of low weight. 

\bibliographystyle{plain}

\end{document}